\documentclass[a4paper,11 pt]{amsart}
\usepackage{srollens-en}[2006/10/02]

\title[Nilmanifolds]{Lie-algebra Dolbeault-cohomology and small deformations of nilmanifolds}

\author{S\"onke Rollenske}

\address{Dr. S\"onke Rollenske\\Department of Mathematics\\
 Imperial College London \\
 SW7 2AZ London\\
 United Kingdom}
\email{s.rollenske@imperial.ac.uk}

\renewcommand{\lg}{\ensuremath{\gothg}}

\newcommand{\einsnull}[1]{{{#1}^{1,0}}}
\newcommand{\nulleins}[1]{{{#1}^{0,1}}}
\newcommand{\nijen}{\mathcal N}
\newcommand{\defobj}[1]{\textbf{#1}}
\newcommand{\Kur}{\mathrm{Kur}}

\newcommand{\refb}[1]{{\upshape (\ref{#1})}}

\begin{document}

\begin{abstract}
We consider nilmanifolds with left-invariant complex structure and prove that  small deformations of such structures are again left invariant if the Dolbeault-cohomology of the nilmanifold can be calculated using left-invariant forms. By a result of Console and Fino this is generically the case. Our main tool is an analog of Dolbeault-cohomology for Lie-algebras with complex structure.

\end{abstract}

\maketitle

\tableofcontents
\section*{Introduction}

Let $X=(M,J)$ a compact complex manifold  regarded as a differentiable manifold together with an integrable almost complex structure. Then it is natural to ask which other complex structures can be obtained by deforming the given complex structure $J$. In general this question is difficult to answer but at least for small deformations there is a universal theory due to Kodaira, Spencer and Kuranishi \cite{kod-sp58, kuranishi62}.

Nevertheless, one rarely has a detailed understanding of all small deformations of a given manifold. One class of manifolds where one can hope for such are so-called nilmanifolds with left-invariant complex structure, i.e.,  compact quotients of real nilpotent Lie-groups  equipped with a left-invariant complex structure.

Nilmanifolds with left-invariant complex structure provide an important source for examples in complex
differential geometry. Among these are the so-called Kodaira-Thurston manifolds, historically the first examples known to  admit both a complex structure and a symplectic structure but no K\"ahler structure. In fact, a nilmanifold $M$ admits a K\"ahler structure if and only if 
it is a complex torus \cite{ben-gor88,hasegawa89}.

Nilmanifolds will be described by a  triple $(\lg, J, \Gamma)$ where $\lg$ is the nilpotent Lie-algebra associated to a simply connected nilpotent Lie-group $G$, $J$ is an integrable complex structure on $\lg$ (see Section \ref{basicdefin}) and $\Gamma\subset G$ is a (cocompact) lattice. The datum of either $\lg$ or $\Gamma$ (considered as an abstract group) determines $G$ up to unique isomorphism.

The general philosophy is that the geometry of the compact, complex manifold $M_J=(\Gamma\backslash G,J)$ should be completely determined by  the linear algebra of $\lg$, $J$ and the $\IQ$-subalgebra generated by $\log \Gamma\subset \lg$. 

In order to control small deformations using Kuranishi's  theory we have to get a good grip on the Dolbeault-cohomology of the holomorphic  tangent bundle.
In Section \ref{LDC} we set up a Lie-algebra Dolbeault-cohomology with values in integrable modules (Definition \ref{definintegrable}) and prove an analogue of Serre-Duality in this context. Since it is known that for nilmanifolds the usual Dolbeault-cohomology $H^{p,q}(M)=H^q(M, \Omega^p_M)$ can (nearly always) be calculated using invariant forms \cite{con-fin01, cfgu00},  we can identify the cohomology of the tangent bundle with the the cohomology of the complex
\[ 0\to \einsnull\lg \overset{\delbar}{\to} \nulleins{\lg^*}\tensor \einsnull\lg\overset{\delbar}{\to} \Lambda^2\nulleins{\lg^*}\tensor \einsnull\lg\overset{\delbar}{\to} \dots \]
as explained in Section \ref{liedolbeault}.
The explicit description of the cohomology of the tangent bundle will then be used in Section \ref{smalldefo} to prove our main result.
\begin{custom}[Theorem \ref{invariantdeformation}]
Let $M_J$ $(\lg, J, \Gamma)$ be a nilmanifold with left-invariant complex structure such that 
\begin{equation}\label{cohomcond}
\iota:H^{1,q}((\lg,J),\IC)\to H^{1,q}(M_J) \text{ is an isomorphism for all $q$}.\tag{$\ast$}
\end{equation}
Then all small deformations of the complex structure $J$ are also left-invariant complex structures. More precisely, the Kuranishi family contains only left-invariant complex structures. 
\end{custom}
This generalises the analogous result for  abelian complex structures due to Console, Fino and Poon \cite{con-fin-poon06} (see also \cite{mpps06}).
Deformations in the large will be studied in \cite{rollenske08d}.

Since there are no examples known for which \refb{cohomcond} does not hold, it is widely believed that the following question has a positive answer:
\begin{custom}[Question 1]
Does \refb{cohomcond} hold for every left-invariant complex structure on a nilmanifold?
\end{custom}
We hope that Lie-algebra Dolbeault-cohomology theory will turn out to be useful also in other context, for example in the study of the  differential graded algebras arising from nilpotent Lie-algebras with a complex structure (see e.g. \cite{0708.3442v2, math.DG/0610768}).

\subsubsection*{Acknowledgements.} This work is part of my PhD-thesis \cite{rollenske07}. I would like to express my gratitude to my adviser Fabrizio Catanese for suggesting this research, constant encouragement and several useful and critical remarks. Simon Salamon gave  valuable references to the literature and helped to improve the presentation at several points. An invaluable bibliographic hint due to Oliver Goertsches opened a new perspective on the problem. Several suggestions by the referee helped to improve the presentation of the results.

During the revision of the paper I was visiting Imperial College London supported by a DFG Forschungsstipendium.

\section{Lie-algebra Dolbeault-cohomology}\label{LDC}

The aim of this section is to set a Dolbeault-cohomology theory for modules over Lie-algebras with complex structure and prove Serre-duality in this context. Our main application is the calculation of the cohomology groups of the tangent sheaf for nilmanifolds with a left-invariant complex structure in Section \ref{smalldefo} but perhaps the notions introduced are of independent interest. 

After recalling the basic definitions we define the notion of (anti-)integrable module in Section \ref{intreps} and derive the elementary properties which are interpreted in geometric terms in Section \ref{mod-VB}. Section \ref{liedolbeault} is again devoted to algebra when we set up our machinery of Lie-algebra Dolbeault-cohomology. The geometric implications of this theory will be explained in Section \ref{invariantcohomology}.

\subsection{Lie-algebras with a complex structure}\label{basicdefin}
Let $\lg$ be a finite dimensional real Lie-algebra and $J$ an almost complex structure on the underlying real vector space, i.e., $J$ is an endomorphism of $\lg$ such that $J^2=-Id_\lg$; setting $ix=Jx$ for $x\in \lg$ this makes $\lg$ into a complex vector space.
We can decompose the complexified Lie-algebra $\lg_\IC={\gothg}^{1,0}\oplus \gothg^{0,1}$ into the $\pm  i$-eigenspaces of the $\IC$-linear extension of $J$ and every decomposition $\lg_\IC=U\oplus \bar U$ gives rise to a unique almost complex structure $J$ such that  $\einsnull{\lg}=U$.
 
Usually we will use small letters $x,y,\dots$ for elements of $\lg$ and capital letters $X,Y,\dots$ for elements in $\einsnull\lg$. Elements in $\nulleins\lg$ will be denoted by $\bar X,\bar Y,\dots$ using complex conjugation; we will use the same symbol for a linear map and its complexification.

The exterior algebra of the dual vector space $\lg^*$ decomposes as
\[\Lambda^k\lg^*_\IC=\bigoplus_{p+q=k}\Lambda^p{\gothg^*}^{1,0}\tensor \Lambda^q{\gothg^*}^{0,1}=\bigoplus_{p+q=k}\Lambda^{p,q}{\gothg^*}\]
and we have $\overline{\Lambda^{p,q}{\gothg^*}}=\Lambda^{q,p}{\gothg^*}$.
A general reference for the linear algebra coming with a complex structure is  \cite{Huybrechts} (Section 1.2).

\begin{defin}\label{defintegrable}
An almost complex structure $J$ on a real Lie-algebra $\lg$ is said to be \defobj{integrable} if the Nijenhuis condition
\begin{equation}\label{nijenhuis}
 [x,y]-[Jx,Jy]+J[Jx,y]+J[x,Jy]=0
\end{equation}
holds for all $x,y\in \lg$ and in this case we call the pair $(\lg , J) $  a \defobj{Lie-algebra with complex structure}. 
\end{defin}

Hence by a complex structure on a Lie-algebra we will always mean an integrable one. Otherwise we will speak of almost complex structures. We will mainly be concerned with nilpotent Lie-algebras.

\begin{rem}
\begin{enumerate}
\item The real Lie-algebra $\lg$ has the structure of a complex Lie-algebra induced by $J$ if and only if $J[x,y]=[Jx,y]$ holds for all $x,y\in\lg $ and $J$ is then automatically integrable in the above sense.
\item One  shows that $J$ is integrable if and only if $\lg^{1,0}$ is a subalgebra of $\lg_\IC$ with the induced bracket. 
\item If $G$ is a real Lie-group with Lie-algebra $\lg$ then giving a left-invariant almost complex structure on $G$ is equivalent to giving an almost complex structure $J$ on $\lg$ and $J$ is integrable if and only if it is integrable as an almost complex structure on $G$. It then induces a complex structure on $G$ by the Newlander-Nirenberg theorem (\cite{Kob-NumII}, p.145) and $G$ becomes a complex manifold. The elements of $G$ act holomorphically by multiplication on the left but $G$ is not a complex Lie-group in general. 
\end{enumerate}
\end{rem}

\subsection{Integrable representations and modules}\label{intreps}

For the whole section let $(\lg, J) $ be a Lie-algebra with complex structure.
%Several times we will refer to the Nijenhuis tensor \refb{nijenhuis} which was defined in Section \ref{basicdefin}. 

A left $\lg$-module structure on a vector space $E$ is given by a bilinear map 
\[\lg\times E \to E \qquad (x,v)\mapsto x\cdot v\] 
such that $[x,y]\cdot v=x\cdot y\cdot v-y\cdot x\cdot v$. Note that this induces a map $\lg \to \End E$, a representation of $\lg $ on $E$. If we want to stress the Lie-algebra structure on $\End E$ (induced by the structure of associative algebra by setting $[a,b]=ab-ba$) we use the notation $\mathfrak{gl}(E)$. A representation or a left module structure corresponds hence to a Lie-algebra homomorphism $\lg\to \mathfrak{gl}(E)$. 

In the sequel we want to combine these notions with complex structures both on $\lg$ and $E$.
Let $(\lg,J)$ be a Lie-algebra with (integrable) complex structure and $E$ a real vector space with (almost) complex structure $I$.
\begin{defin}\label{definintegrable}
A representation $\rho: (\lg,J) \to \End E$ of $\lg$ on $E$ is said to be  \defobj{integrable} if for all $x\in \lg$ the endomorphism of $E$ given by
\[ \nijen(x):=[I, (\rho\circ J)(x)]+I[\rho(x), I]\]
vanishes identically. In this case we say that $(E,I)$ is an \defobj{integrable $(\lg,J)$-module}. We say that $(E,I) $ is \defobj{anti-integrable} if $(E,-I) $, the complex conjugate module, is an integrable $\lg$-module.
A homomorphism of (anti-) integrable $\lg$-modules is a homomorphism of underlying $\lg$-modules which is $\IC$-linear with respect to the complex structures.
\end{defin}

This definition is modeled on the adjoint representation $ad:\lg \to \End\lg$ which is integrable if and only if 
 the Nijenhuis tensor \refb{nijenhuis} vanishes. This is a special case of the next result.
% \end{rem}

\begin{prop}\label{subrep} 
Let   $(E,I)$ a be vector space with complex structure and $\rho:\lg\to \End E$ a
representation. Then the following are equivalent:
\begin{enumerate}
\item $\rho$ is integrable.
\item For all $X\in \lg^{1,0}$ the map $\rho(X)$ has no component in $\Hom(E^{1,0},E^{0,1})$.
\item  $E^{1,0}$ is an invariant subspace under the action of $\lg^{1,0}$.
\item  $\rho\restr{\lg^{1,0}} $ induces a $\IC$-linear representation on $E^{1,0}$ by
restriction.
\end{enumerate}
\end{prop}
\pf
The restriction of $\rho$ to $\einsnull \lg$ is $\IC$-linear by definition since it is the complexification of the real
representation restricted to a complex subspace. Therefore condition (\textit{ii}) is  equivalent to (\textit{iii}) and (\textit{iv}).

 It remains to prove (\textit{i})$\iff$ (\textit{ii}).
 Let $X\in \einsnull{\lg}$ and $V\in \einsnull E$. Using $JX=iX$ and $IV=iV$ we calculate
\begin{align*}
\nijen(X)V&=  ([I, (\rho\circ J)(X)]+I[\rho(X), I])(V)\\
&= (iI\rho(X)-i\rho(X)I+ I\rho(X)I-I^2\rho(X))(V)\\
&= 2iI\rho(X)V+2\rho(X)V
\end{align*}
and see that 
\[\nijen(X)V=0\iff I\rho(X)V=i\rho(X)V\iff\rho(X)V\in\einsnull V.\]
This proves (\textit{i})$\implies$ (\textit{ii}). Vice versa assume that (\textit{ii}) holds. We decompose the elements
in $E_\IC$ respectively in $\lg_\IC$ into their $(1,0)$ and $(0,1)$
parts. By the above calculation and its complex conjugate (the representation and hence the bracket are real and commute with complex conjugation) it remains to consider the \emph{mixed} case. We have for all
$X,V$ as above
\begin{align*}
\nijen(X)\bar V&= (iI\rho(X)-i\rho(X)I+ I\rho(X)I-I^2\rho(X))(\bar
 V)\\
 &= iI\rho(X)\bar V-\rho(X)\bar V - i  I\rho(X)\bar V +\rho(X)\bar
 V\\
 &=0
\end{align*}
and hence $\rho$ is integrable.\qed

\begin{cor} \label{antidelbarrep}Let   $(E,I)$ a be vector space with complex structure and $\rho:\lg\to \End E$ a
representation. Then the following are equivalent:
\begin{enumerate}
\item $\rho$ is anti-integrable.
\item For all $\bar X\in \lg^{0,1}$ the map $\rho(\bar X)$ has no component in $\Hom(E^{1,0},E^{0,1})$.
\item  $E^{1,0}$ is an invariant subspace under the action of $\lg^{0,1}$.
\item  $\rho\restr{\lg^{0,1}} $ induces a $\IC$-linear representation on $E^{1,0}$ by
restriction.
\end{enumerate}
\end{cor}
\pf Exchange $I$ by $-I$ in the proof of Proposition \ref{subrep}.\qed

\begin{prop} \label{delbarrep} 
Let  $\rho$ be an integrable representation on $(E,I)$. The bilinear map $\delta$ given by 
\[ \delta:\nulleins\lg\times \einsnull E\overset{\rho}{\to} E_\IC\overset{pr}{\to} \einsnull E\quad (\bar X, V)\mapsto \einsnull{(\rho(\bar
X)V)}\]
induces a $\IC$-linear  representation of $\nulleins \lg$ on $\einsnull E$.
\end{prop}
\pf
The map is clearly complex bilinear and it remains to prove the
compatibility with the bracket. Let $\bar X, \bar Y\in \nulleins \lg $
and $V\in \einsnull E$ be arbitrary. Note that $\rho(\bar Y)V= \delta (\bar Y,V)+\nulleins{(\rho(\bar Y)V)}$. Then
\begin{align*}
&\delta([\bar X, \bar Y],V)=\einsnull{(\rho([\bar X, \bar Y])V)}\\
=& \einsnull{\left(\rho(\bar X)\rho(\bar Y)V-\rho(\bar Y)\rho(\bar X)V\right)}\\
=& \einsnull{\left(\rho(\bar X)\left(\delta(\bar Y,V)+\nulleins{(\rho(\bar Y)V)}\right)
 -\rho(\bar Y)\left(\delta(\bar X,V)+\nulleins{(\rho(\bar X)V)}\right)\right)}\\
=&\einsnull{\left(\rho(\bar X)\delta(\bar Y,V) -\rho(\bar Y)\delta(\bar
  X,V)\right)}+\einsnull{(\underbrace{\rho(\bar X)\nulleins{(\rho(\bar Y)V)}-\rho(\bar Y)\nulleins{(\rho(\bar X)V)}}_\text{of type (0,1)})}\\
=&\delta(\bar X,\delta (\bar Y, V))-\delta(\bar Y,\delta (\bar X,V)).
\end{align*}
Here we used that the action of $\nulleins \lg$ maps $\nulleins E$ to $\nulleins E$ which is the complex conjugate of Proposition \ref{subrep} (iii). Hence $\delta$ induces a $\nulleins \lg $-module structure on
$\einsnull E$ as claimed.\qed

We want to clarify the relation between integrable and anti-integrable modules.
\begin{lem}
Let $(E,I)$ be an integrable  $(\lg, J)$-module. Then the dual module  with the induced $\lg$-module structure is anti-integrable. 
\end{lem}
\pf  If $x\in \lg$ and $\phi\in E^*$ then the induced module structure is given by  $(x\cdot\phi)(v)=-\phi(xv)$  for $v\in E$. We have to show that for $\bar X\in \nulleins \lg$ and $\Phi\in\einsnull{ E^*}$ the map $(\bar X\cdot\Phi)$ annihilates $\nulleins E$. But if $\bar V$ is in $\nulleins E$ then by the above proposition $\bar X\bar V\in \nulleins E$ and $\Phi(\bar X\bar V)=0$.\qed

\begin{rem}\label{badcategory}
The above result seems unnatural only at first sight. If we consider $E$ as left $\mathfrak{gl} (E)$-module in the canonical way then the complex structure $J\in \End E$ acts on the left. The dual vector space $E^*$ comes with a natural action of $\mathfrak{gl}( E)$ on the right:
\[ \phi \cdot A:= \phi\circ A\qquad\text{ for }A\in \mathfrak{gl}( E), \phi \in E^*\]
and the complex structure of $E^*$ is given exactly in this way  $I^*\phi=\phi\circ I$.

In order to make $E^* $ a left $\mathfrak{gl}(E)$-module we have to change sign $A\cdot \phi := -\phi\circ A$; changing the sign of the complex structure corresponds to complex conjugation.

Integrable modules do not behave well under standard linear algebra operation like the tensor product. The reason is simply that we have to work over $\IC$ if we want to keep the complex structure on our vector spaces and over $\IR$ if we want to preserve the $\lg$ action, since this action is not $\IC$-linear in general. 
\end{rem}

\subsection{Integrable modules and vector bundles}\label{mod-VB}
Now we want to relate the notion of integrable $\lg$-module to geometry. First we forget the complex structures and look at the differentiable situation:

Let $\lg$ be a real Lie-algebra and $G$ the corresponding simply connected Lie-group. Let $E$ be a (left) $\lg$-module and $\Gamma\subset G$ a co-compact discrete subgroup. Then $E \times G$ is the trivial bundle with an  action of $G$ on the left, given by the representation of $\lg$ on $E$. 

If we take the quotient by the action of the subgroup $\Gamma$ then the result is a homogeneous, flat vector bundle on $M=\Gamma\backslash G$.

Another possibility to look at this situation is the following: the representation of $\lg$ on $E$ gives rise to a central extension, the semi-direct product $E\rtimes \lg$. The vector space underlying this Lie-algebra is $E\oplus \lg $ with the Lie-algebra structure given by $[(v,x), (w, y)]= (x\cdot w-y\cdot v, [x,y])$ for $(v,x), (w, y)\in E\oplus \lg$.

 Regarding the real vector space $E$ as a commutative Lie-group, we an the exact sequence of Lie-groups
\[ 0 \to E\to E\rtimes G\to G\to 1.\]

Now we take the complex structures $(\lg,J)$ and $(E,I)$ into account: using left-translations $J$ induces a left-invariant integrable almost complex structure on $G$ and the quotient $M_J:=(\Gamma\backslash G, J)$  is a compact complex manifold on which the normaliser of $\Gamma$ acts holomorphically on the left and the whole group acts differentiably on the right.

%the complex structure on $\lg$ induces a left-invariant almost complex structure on $G$
% by declaring $\einsnull{(T_gG)}:= {l_g}_* \einsnull \lg $. Here we identify $\lg $ with the tangent space at the identity and $l_g$ is the left multiplication by an element $g\in G$.
% 
% This almost complex structure is integrable if and only if  $J$ is integrable in the sense of \refb{nijenhuis}
%  and we assume this to be the case in what follows; thus $(G,J)$ is	 a complex manifold.
% 
% By definition  left multiplication $l_g$  is holomorphic but in general the multiplication on the right will not be holomorphic. If this is indeed the case, then $J$ is bi-invariant and $G$ is  a complex Lie-group. 

Given a $(\lg, J)$-module with a complex structure $(E,I)$ we would like to define the structure of a complex vector bundle on the differentiably trivial vector bundle $E\times G$, that descends to a flat homogeneous holomorphic vector bundle on $M_J$. This is possible if $(E,I)$ is integrable:
\begin{lem}
The $(\lg, J)$-module $(E,I)$ is integrable if and only if $I\times J$ induces a left-invariant complex structure on the Lie-group $E\rtimes G$. 
\end{lem}
\pf
 We have to calculate the Nijenhuis tensor \refb{nijenhuis} for $K=I\times J$ on the Lie-algebra $E\rtimes \lg$. Let 
$(v,x), (w, y)$ be in $E\rtimes \lg$. Then
\begin{align*}
& [(v,x), (w, y)]-[K(v,x), K(w, y)]+K[K(v,x), (w, y)]+K[(v,x), K(w, y)]\\
=& (x\cdot w-y\cdot v, [x,y])- (Jx\cdot Iw-Jy\cdot Iv, [Jx,Jy]) \\
&+(I(Jx\cdot w-y\cdot Iv), J[Jx,y]) + (I(x\cdot Iw-Jy\cdot v), J[x,Jy]).
\end{align*}
The second component yields the Nijenhuis tensor for $J$ and hence vanishes since we assumed $J$ to be integrable.

Using the notation $\rho(x)$ for the element in $\End E$ corresponding to the action of $x$ on $E$ the we expand $\kn(x)$ as
\begin{align*}
\kn(x)=&I[\rho(x), I]+[I, \rho(Jx)]\\
=&-{I}^2\rho(x) +I\circ\rho(x)\circ I+[I, \rho(Jx)]\\
=&\rho(x)-\rho(Jx)\circ J+I\circ\rho(Jx)+I\circ\rho(x)\circ I,
\end{align*}
and so  the first component becomes
\begin{align*}
 &x\cdot w-y\cdot v- Jx\cdot Iw+Jy\cdot Iv \\
&\qquad\qquad+I(Jx\cdot w-y\cdot Iv)+(I(x\cdot Iw-Jy\cdot v)\\
=&(\rho(x)-\rho(Jx)\circ J+I\circ\rho(Jx)+I\circ\rho(x)\circ I)w\\
&\qquad\qquad-(\rho(y)-\rho(Jy)\circ J+I\circ\rho(Jy)+I\circ\rho(y)\circ I)v\\
=&\kn(x)w-\kn(y)v.
\end{align*}
  Setting $v=0$ we see that $K$ is integrable if and only if $(E,I)$ is an integrable $(\lg, J)$-module.\qed

\begin{rem}
\begin{enumerate}
\item Note that the left-invariance of a complex structure is not preserved by standard   vector bundle operations for essentially the same reason as explained in Remark \ref{badcategory}.
\item Even if a vector bundle $E$ is equipped with an integrable left-invariant structure left-invariant sections are not necessarily holomorphic. In fact, the action of each $g\in G$ on $E$ is holomorphic but there is no natural way to speak about \emph{the action varying holomorphically} since $G$ is not a complex Lie-group in general.
\end{enumerate}
\end{rem}

\subsection{Lie-algebra Dolbeault-cohomology}\label{liedolbeault}

In this section we want to define a cohomology theory for Lie-algebras with a complex structure with values in a finite-dimensional integrable module. In the notation we will often suppress the complex structures.

Recall that the cohomology groups of a Lie-algebra $\lg$ with values in a $\lg$-module $E$ are defined as the right derived functors of the  functor of invariants (\cite{Weibel}, Chapter 7)
\[E\mapsto E^\lg=\{m\in E\mid x\cdot m=m  \text{ for all } x\in \lg\}.\]

The cohomology groups can be calculated using the so-called Chevalley complex for $E$:
\[ 0\to E\overset{d_0}{\to} \lg^*\tensor E\overset{d_1}{\to} \Lambda^2 \lg^* \tensor E\overset{d_2}{\to} \dots {\to} \Lambda^{dim \lg} \lg^*\tensor E\to 0\]
with differential given by
\begin{multline}\label{ch-diff}
(d_k\alpha)(x_1, \dots , x_{k+1}):=\sum_{i=1}^{k+1} (-1)^{i+1} x_i (\alpha(x_1, \dots ,\hat x_i, \dots , x_{k+1}))\\
+\sum_{1\leq i <j\leq k+1} (-1)^{i+j} \alpha([x_i,x_j], x_1, \dots, \hat x_i, \dots, \hat x_j,\dots , x_{k+1}).
\end{multline}
It was originally introduced by Chevalley and Eilenberg \cite{Knapp,che-eil48} and its elements  can be interpreted as left-invariant differential forms in the the geometric context.

Now let $(\lg,J)$ be a Lie-algebra with complex structure and let $(E,I)$ be a  finite dimensional, integrable (resp. anti-integrable)  $\lg$-module. By Proposition \ref{delbarrep} (resp. Corollary \ref{antidelbarrep}) we have a representation of $\nulleins\lg $ on $\einsnull E$ given by
$(\bar X, V)\mapsto \einsnull{(\rho(\bar X)V)}$ (resp. $(\bar X, V)\mapsto \rho(\bar X)V$). Together with the representation of $\nulleins \lg$ on $\Lambda^p\einsnull{ \lg^*}$ induced by the adjoint representation we obtain a $\nulleins\lg$-module structure on  $\Lambda^p\einsnull{\lg^*}\tensor\einsnull E$. 

\begin{defin}
Let $(\lg,J)$ be a Lie-algebra with complex structure and let $(E,I)$ be a  finite dimensional, integrable (anti-integrable)  $\lg$-module. Then we define
\begin{gather*}
H^{p,q}(\lg, E)=H^{p,q}_{\delbar}((\lg,J),( E,I)):= H^q(\nulleins{\lg}, \Lambda^p\einsnull{\lg^*}\tensor\einsnull E)\\
\end{gather*}
%where $H^k(\nulleins{\lg}, \einsnull E)$ is the Lie-algebra cohomology of $\nulleins\lg$ with values in the $\einsnull E$  with the module structure as in Proposition \ref{delbarrep} (as in Corollary \ref{antidelbarrep}).
We call $ H^k_{\delbar}(\lg, E):=H^{0,k}(\lg, E)$ the $k$-th \defobj{Dolbeault-cohomology group of $\lg$ with values in $E$}.
\end{defin}

\begin{exam}\label{trivialmodule}
Consider for $(\lg, J)$ as above $\IC$ as the trivial $\lg_\IC$-module. Then the associated Chevalley differential on the exterior algebra $\Lambda^\bullet \lg^*_\IC$ decomposes into $d=\del+\delbar$ since $J$ is integrable and we can consider the double complex $(\Lambda^{p,q}\lg^*, \delbar, \del)$. 

The adjoint action of $\lg$ on itself yields  an anti-integrable $\lg$-module structure on $\lg^*$, hence a $\nulleins\lg$-module structure on  $\Lambda^p\einsnull{\lg^*}$. It is now easy to see that the columns of the above double complex 
\[0\to \Lambda^{p,0}\lg^*\to\Lambda^{p,1}\lg^*\to \dots\]
 are the Chevalley complexes calculating $H^q(\nulleins\lg, \Lambda^{p,0}\lg^*)= H^{p,q}_{\delbar} (\lg, \IC)$.
\end{exam}

Now we want to develop some kind of Hodge theory for our Dolbeault-cohomology which we model on the usual Hodge theory for holomorphic vector bundles as it can be found for example in the book of Huybrechts \cite{Huybrechts}.
Let $2n$ be the real dimension of $\lg$.

First of all we choose an euclidean structure $g=\langle-,-\rangle$ on the real vector space underlying $\lg$ which is compatible with the given complex structure $J$ in the sense that $\langle J-,J-\rangle=\langle-,-\rangle$. Let $vol$ be the volume form, i.e., the unique element in $\Lambda^{2n} \lg^*$ inducing the same orientation as $J$ and of length one in the induced metric on $\Lambda^\bullet \lg^*$. We define the Hodge-$\ast$-operator, which is an isometry on $\Lambda^\bullet\lg^*$,  by 
\[ \alpha\wedge \ast\beta = \langle\alpha, \beta\rangle vol \quad \text{for all } \alpha, \beta \in \Lambda^\bullet\lg^*.\]
On the complexified vector space $\lg_\IC$ we have a natural induced hermitian product $\langle-,-\rangle_\IC$ and a $\ast$-operator determined by
\[ \alpha\wedge \ast\bar\beta = \langle\alpha, \beta\rangle_\IC vol \quad \text{for all } \alpha, \beta \in \Lambda^\bullet\lg^*_\IC.\]
which maps $(p,q)$-forms to $(n-p, n-q)$-forms.

We want now to define a star operator also on $\Lambda^{\bullet,\bullet}\lg^*\tensor \einsnull E$. For this purpose we choose an euclidean product on  $E$ compatible with the complex structure $I$ which induces a hermitian structure $h$ on $\einsnull E$. We consider $h$ as an $\IC$-antilinear isomorphism $h:\einsnull E\isom \einsnull{E^*}$.  Then
\[ \bar\ast_{ E} :\Lambda^{p,q}\lg^*\tensor \einsnull E\to \Lambda^{n-p,n-q} \tensor\einsnull  E^*\]
is defined by  $\bar\ast_{ E}(\alpha\tensor s)= \overline{\ast\alpha}\tensor h(s)=\ast(\bar\alpha)\tensor h(s)$.
Let $(-,-)$ be the hermitian form on $\Lambda^{\bullet, \bullet}\lg^*\tensor \einsnull E$ induced by $g$ and $h$. Then $\bar\ast_{ E}$ is a $\IC$-antilinear isomorphism depending on our choice of $g$ and $h$ and the identity
\[ (\alpha, \beta)vol = \alpha\wedge\bar\ast_E\beta\]
holds for $\alpha, \beta\in \Lambda^{p,q}\lg\tensor \einsnull E$, where $"\wedge"$ is the exterior product for the elements in $\Lambda^{\bullet, \bullet}\lg$ and the evaluation map $\einsnull E\tensor \einsnull E^*\to \IC$ on the module part.
It is not difficult to verify  that one has $\bar\ast_{ E^*}\circ\bar\ast_{ E}=(-1)^{p+q}$ on $\Lambda^{p,q}\lg\tensor \einsnull E$.

\begin{defin}\label{laplacedefin}
Let $(E,I)$  be an (anti-)integrable $(\lg, J)$-module. The operator $\delbar_E^*:\Lambda^{p,q}\lg\tensor \einsnull E\to\Lambda^{p,q-1}\lg\tensor \einsnull E$ is defined as
\[\delbar_E^*:= -\bar\ast_{E^*}\circ \delbar_{E^*} \circ \bar\ast_E.\]
Let  $\Delta_E:=\delbar_E^*\delbar_E+\delbar_E\delbar_E^*$ be the \defobj{Laplace operator} on $\Lambda^{p,q}\lg\tensor \einsnull E$. We call an element $\alpha$ harmonic if $\Delta_E(\alpha)=0$ and denote by $\kh^{p,q}(\lg,E)$ the space of harmonic elements (where we omit $g$ and $h$ from the notation).
 \end{defin}
 
Observe that $\bar\ast_E$ induces a $\IC$-antilinear isomorphism 
\[\bar\ast_E : \kh^{p,q}(\lg,E)\isom \kh^{n-p,n-q}(\lg,E^*).\]

\begin{prop} If $H^{2n}(\lg_\IC, \IC)=\IC$, where $\IC$ is considered as the trivial $\lg$-module, then the operator $\delbar_E^*$ is adjoint to $\delbar_E$ with respect to the metric induced by $g$ and $h$. In this case $\Delta_E$ is self-adjoint.
\end{prop}
The condition on the cohomology is somehow the equivalent of Stokes theorem as will be seen in the proof.

\pf The second assertion is a consequence of the first one which in turn is proved by the following calculation:

First of all note that the assumption $\Lambda^{2n}\lg^*_\IC\isom \IC$ implies that $d_{2n-1}=0$ in $\Lambda^\bullet \lg^* $, the Chevalley complex of the trivial module.  Hence the same holds for $\delbar : \Lambda^{n, n-1}\lg^*\to \Lambda^{n,n}\lg^*$. For $\alpha \in \Lambda^{p,q}\lg^*\tensor \einsnull E$ and $\beta \in \Lambda^{p,q+1}\lg^*\tensor \einsnull E$ we have
\begin{align*}
(\alpha, \delbar^*_E\beta)vol &= -(\alpha, \bar\ast_{E^*}\circ \delbar_{E^*} \circ \bar\ast_E\beta)vol\\
&= -\alpha\wedge \bar\ast_E\bar\ast_{E^*}\delbar_{E^*}\bar\ast_E\beta\\
&= (-1)^{n-p+n-q-1} \alpha\wedge\delbar_{E^*}\bar\ast_E\beta\\
&= -\delbar(\alpha\wedge\bar\ast_E\beta)+\delbar_E\alpha\wedge\bar\ast_E\beta\\
&= \delbar_E\alpha\wedge\bar\ast_E\beta\\
&= (\delbar_E\alpha, \beta)vol.
\end{align*}
Here we used the identity 
\[\delbar(\alpha\wedge\bar\ast_E\beta)=\delbar_E\alpha\wedge\bar\ast_E\beta+(-1)^{p+q} \alpha\wedge\delbar_{E^*}\bar\ast_E\beta\]
that follows form the Leibniz rule in the exterior algebra and the fact that the evaluation map $\einsnull E\tensor \einsnull E^*\to \IC$ is a map of $\nulleins \lg $-modules. \qed
\begin{rem}
We have always $H^{2n}(\lg_\IC , \Lambda^{2n}\lg_\IC)=\IC$ (See \cite{Weibel}, Exercise 7.7.2). Hence the assumptions of the theorem hold if $\lg $ acts trivially on $\Lambda^{2n}\lg$, in particular if $\lg $ is nilpotent.
\end{rem}
Here are some standard consequences:
\begin{cor}
If $H^{2n}(\lg_\IC, \IC)=\IC$ then an element $\alpha \in \Lambda^{p,q}\lg\tensor E$  is harmonic if and only if $\alpha$ is $\delbar_E$ and $\delbar_E^*$ closed.
\end{cor}
\pf Standard argument.\qed
\begin{cor}[Hodge-decomposition]
Let $(E,J)$ be an (anti-)integrable module  over the Lie-algebra with complex structure $(\lg,J)$ both equipped with  compatible euclidean products. If $H^{2n}(\lg_\IC, \IC)=\IC$ then there is a orthogonal decomposition 
\[\Lambda^{p,q}\tensor \einsnull E =\delbar_E( \Lambda^{p,q-1}\tensor \einsnull E)\oplus \kh^{p,q}(\lg,E)\oplus 
\delbar_E^*(\Lambda^{p,q+1}\tensor \einsnull E).\]
\end{cor}
\pf Since everthing is finite dimensional this follows trivially from the above.\qed
\begin{cor}
If $H^{2n}(\lg_\IC, \IC)=\IC$ then the natural projection
\[\kh^{p,q}(\lg,E)\to H^{p,q}(\lg,E)\]
 is bijective.
\end{cor}
\begin{theo}[Serre-Duality]
Let $(\lg, J)$ be a Lie-algebra with complex structure such that $H^{2n}(\lg_\IC, \IC)=\IC$ and $(E,I)$ an (anti-)integrable $\lg$-module. Then the paring
\[H^{p,q}(\lg, E)\times H^{n-p,n-q}(\lg, E^*)\to \IC\cdot vol\isom \IC \quad (\alpha, \beta)\mapsto \alpha \wedge \beta\]
 is well defined and non degenerate.% where $"\wedge"$ is as above.
\end{theo}
\pf
Fix hermitian structures on $E$ and $\lg$ respectively. Then consider the pairing 
\[\kh^{p,q}(\lg, E)\times \kh^{n-p,n-q}(\lg, E^*)\to \IC\cdot vol\isom \IC.\]

We claim that for any non-zero $ \alpha\in \kh^{p,q}(\lg, E)$ there exists an element $\beta \in \kh^{n-p, n-q}(\lg, E^*)$ such that $\alpha \wedge \beta \neq 0$. Indeed, choosing $\beta= \bar\ast_E\alpha$ we have 
\[\alpha\wedge\beta=\alpha\wedge\bar\ast_E\alpha=(\alpha, \alpha)vol =\|\alpha\|^2 vol\neq 0.\]
This proves that the pairing in non degenerate.\qed

\begin{cor}\label{serreduality} Let $(\lg, J)$ be a Lie-algebra with complex structure such that $H^{2n}(\lg_\IC, \IC)=\IC$
For any (anti-)integrable $(\lg, J)$-module there exist natural isomorphisms
\begin{gather*}
H^{p,q}(\lg, E)\isom H^{n-p, n-q}(\lg, E^*)^*\\
\intertext{and if $\Lambda^n\lg^*$ is the trivial $\lg$-module}
H^q_{\delbar}(\lg, E)\isom H^{n-q}_{\delbar}(\lg, E^*)^*.
\end{gather*}
\end{cor}

\subsection{Cohomology with invariant forms}\label{invariantcohomology}
 We are now going to translate our results on  Lie-algebra Dolbeault-cohomology to the geometric situation.

Recall the situation considered in Section \ref{mod-VB}: let $(\lg, J)$ be a real Lie-algebra with complex structure of real dimension $2n$ and  $(E,I)$ an integrable $(\lg,J)$-module. Let $G$ be the simply connected Lie-group associated to $\lg$ endowed with the left-invariant complex structure induced by $J$. Let $\Gamma$ a uniform lattice  in $G$ and consider the flat, homogeneous, holomorphic vector bundle $\ke$ on $M=\Gamma \backslash G$ constructed by taking the quotient of $E\times G$ by $\Gamma$ acting on the left.

Let $g$ be a euclidean structure on $\lg$ compatible with the complex structure $J$, such that $M$ has volume one with respect to the associated left-invariant metric on $M$. Choose also a euclidean structure on $E$ compatible with the complex structure $I$.

Let $\pi:G\to M$ be the projection. We say that a smooth section  $s\in \ka^{p,q}(M,\ke)$ is invariant if $\pi^*s$ is invariant under the action of $G$ on the left, i.e., $l_g^*(\pi^*s)=\pi^*s$ for all $g\in G$. This makes sense since $\pi^*\ke=E\times G$ is trivial as a smooth vector bundle and in particular $l_g^*\pi^*\ke=\pi^*\ke$. A smooth section $s$ in the trivial bundle $E\times G$ is the pullback of a section of $\ke$ if and only if it is invariant under the action of $\Gamma$ on the left.

The relation between the usual Dolbeault theory for vector bundles on complex manifolds and our theory developed so far is summarised in the following

\begin{prop}
In the above situation $\Lambda^{p,q}\lg\tensor\einsnull E$ can be identified with the subset of invariant, smooth differential form on $M$ with values in the holomorphic bundle $\ke$. 
If in addition $H^{2n}(\lg_\IC, \IC)=\IC$ then  the following holds:
\begin{enumerate}
\item The differential in the Chevalley complex as given in \refb{ch-diff} coincides with the usual differential restricted to invariant forms with values in $\ke$. In particular, if $E$ is the trivial module the decomposition $d=\del+\delbar$ on $\Lambda^{\bullet, \bullet}\lg^*$ coincides with the usual one on the complex manifold $M$.
\item The Chevalley complex associated to the $\nulleins\lg$-module structure on $\einsnull E$ is the subcomplex of invariant forms contained in the usual Dolbeault resolution of the holomorphic vector bundle $\ke$ by smooth differential forms with values in $\ke$.
\item The Hodge-$\ast$-operator defined on $\Lambda^\bullet \lg^*_\IC$ in Section \ref{liedolbeault} coincides with the usual Hodge-$\ast$-operator on the exterior algebra of smooth differential forms. The same holds true for the operator $\bar\ast_E$.
\item The operators $\delbar^*_E$ and $\Delta_E$ in Definition \ref{laplacedefin} are the restrictions of the corresponding operators on smooth differential forms. In particular we have an inclusion
\[\kh^{n-p,n-q}(\lg, E)\subset \kh^{n-p,n-q}(M, \ke)\]
where $\kh^{n-p,n-q}(M, \ke)$ are the harmonic $(p,q)$-forms with values in $\ke$ with respect to the chosen left-invariant hermitian structures.
\end{enumerate}
\end{prop}
\pf The first claim is  clear by construction. The Lie bracket on $\lg$ is clearly the restriction of the usual Lie bracket on vector fields on $M$ and also the definition of the differential in \refb{ch-diff} coincides with the usual one for smooth differential forms (see e.g. \cite{Huybrechts}, p. 283). Since $\pi^*\ke$ is differentiably a trivial bundle the same holds for differential forms on $G$ with values in $\pi^*\ke$ and therefore also for sections of $\ke$ itself since we can check this locally. This proves (i) and (ii) using the identification of $\ke$ with $\einsnull E$.

Our reference for the Hodge theory of holomorphic vector bundles is \cite{Huybrechts} (ch. 4.1). Now, recall that we defined our operator $\bar \ast_E$ by the relation 
\[ \alpha\wedge\bar\ast_E\beta = (\alpha, \beta)vol =(\alpha, \beta)\ast 1\]
for $\alpha, \beta\in \Lambda^{p,q}\lg\tensor \einsnull E$ which conicides with the definition for differential forms in $\ka^{p,q}(M,\ke)$ if we consider $\alpha$ and $\beta$ as invariant differential forms on $M$:

The hermitian metric on $\ka^{p,q}(M,\ke)$ is  defined by $ (\alpha, \beta)_M=\int_M(\alpha, \beta)vol $ but if the forms are invariant we have
\[(\alpha, \beta)_M=\int_M(\alpha, \beta)vol=(\alpha, \beta)\int_M vol =(\alpha,\beta)\]
 since we  chose the invariant metric such that the volume of $M$ is one.

Therefore also $\bar\ast_\ke=\bar\ast_E$ and this concludes the proof since the Laplace operator can be described in terms of $\bar\ast_\ke$ and $\delbar$.\qed

\begin{cor}\label{isoduality}
In the above situation we have an inclusion
\[\iota_E:H^{p,q}(\lg,E)\to H^{p,q}(M,\ke)\]
induced by the inclusion on the level of harmonic differential forms. In particular if $\iota_E$ is an isomorphism then so is $\iota_{E^*}:H^{n-p,n-q}(\lg,E^*)\to H^{n-p,n-q}(M,\ke^*)$.
\end{cor}
\pf The first claim is an immediate consequence of (\textit{vi}) in the proposition while the second then follows for dimension reasons from Serre-Duality both on $M$ and for Lie-algebra Dolbeault-cohomology (Corollary \ref{serreduality}).\qed 

We will apply this to the cohomology of nilmanifolds in the next section in order to study the space of infinitesimal deformations.

\section{Nilmanifolds and their small deformations}\label{smalldefo}

The aim of this section is to prove the main result of this paper namely that small deformations of nilmanifolds with left-invariant complex structure carry a left-invariant complex structure under very mild conditions. 

%We will always regard a complex manifolds as a differentiable manifolds together with an integrable almost complex structure.

\subsection{Nilmanifolds with left-invariant complex structure }\label{set-up}

In this section we will present those aspects of the theory of nilmanifolds with left-invariant complex structure which we need to formulate and prove our theorem. A more detailed study of their geometry can be found in  \cite{rollenske08d}.

In the following let $(\lg, J)$ be a nilpotent Lie-algebra with complex structure as in Definition \ref{defintegrable} and $G$ an associated simply-connected Lie-group. By left-translation $J$ defines an integrable almost complex structure on $G$.

Now assume that $\Gamma\subset G$ is a lattice, a discrete cocompact subgroup. A nilpotent Lie-group $G$ contains such a subgroup if and only if its Lie-algebra can be defined over $\IQ$ (see e.g. \cite{Cor-Green}, p. 204). Then the complex structure on $G$ descends to a complex structure on the compact manifolds $M:=\Gamma\backslash G$.

%By a torus we will always mean a compact torus.

\begin{defin}
A compact complex manifold $M_J:=(M,J)$ is called \defobj{nilmanifold with left-invariant complex structure} if there is a nilpotent Lie algebra with complex structure $(\lg,J)$ and a lattice $\Gamma$ in an associated simply-connected Lie-group such that $M_J\isom(\Gamma\backslash G, J)$.
\end{defin}

Since $\Gamma=\pi_1(M_J)$ determines $G$ and hence $\lg$ up to isomorphism (\cite{VinGorbShvart}, p.45, Corollary 2.6),  by abuse of notation we will always identify  $M_J$ with $(\Gamma\backslash G, J)$ and call it a nilmanifolds with left-invariant complex structure of type $(\lg, \Gamma)$.

\begin{rem}
 \begin{enumerate}
\item 
A nilmanifold $M_J$ of type $(\lg, \Gamma)$ with left-invariant complex structure is K\"ah\-leri\-an if and only if $\lg$ is abelian and $M$ is a complex torus \cite{ben-gor88,hasegawa89}.  It can  be arbitrarily far from being  K\"ahler in the sense that the Fr\"olicher spectral sequence may be arbitrarily non-degenerate \cite{rollenske07a}.
\item   If $[\einsnull{\lg}, \nulleins{\lg}]=0$ then $J$ is called complex parallelisable (bi-invariant), $(G,J)$ is a complex Lie group and $M_J$ is complex parallelisable. They have interesting arithmetic properties \cite{winkelmann98}. The other extremal case are abelian complex structures defined by  $[\einsnull{\lg}, \einsnull{\lg}]=0$ studied for example in \cite{mpps06, con-fin-poon06}.
\item All (iterated) principal holomorphic torus bundles can be described as nilmanifolds with left-invariant complex structure but not all of them admit such a geometric description.
\item 
For nilpotent Lie-groups the exponential map $\exp: \lg \to G$ is a diffeomorphism  and the preimage of $\Gamma$ in $\lg$ generates a rational subalgebra $\lg_\IQ\subset \lg$ (\cite{Cor-Green}, p.204). The complex structure $J$ is called $\Gamma$-rational if it maps $\lg_\IQ$ to itself; in this case one has more control over the geometry of $M_J$ \cite{con-fin01}.
\end{enumerate}
\end{rem}

\subsection{A parameter-space for left-invariant complex structures}

We want to parametrise the space of left-invariant complex structures on a given nilpotent Lie-algebra $\lg$. Let $2n$ be the real dimension of $\lg$. 

A complex structure $J$ is uniquely determined by specifying the $(0,1)$-subspace $\bar V\subset\lg_\IC$ and the integrability condition can be expressed as $[\bar V,\bar V]\subset \bar V$. Hence we write (like in \cite{salamon01})
\[\kc(\lg):= \{\bar V\in \IG(n , \lg_\IC)\mid V\cap \bar V=0, [\bar V,\bar V]\subset \bar V\}\]
where $\IG(n , \lg_\IC)$ is the Grassmannian of $n$-dimensional complex subspaces of $\lg_\IC$.
Recall that its tangent space at a point $\bar V$ is \[T_{\bar V} \IG(n , \lg_\IC)=\Hom_\IC(\bar V, \lg_\IC/\bar V)\isom {\bar V}^*\tensor V\isom \nulleins \lg^*\tensor \einsnull \lg\]
if we endow $\lg$ with the complex structure $J_{\bar V}$ induced by $\bar V$.

In general it is a difficult question to decide if $\kc(\lg)$ is non-empty for a given Lie-algebra $\lg$. For the next paragraph we will assume this to be the case.

Now fix a simply connected nilpotent Lie-group $G$ with Lie-algebra $\lg$. We want to describe a family of complex manifolds $\pi:\km(\lg)\to \kc(\lg)$ such that over every point $\bar V\in \kc(\lg)$ the fibre $\inverse\pi(\bar V)$ is the manifold $G$ with the left-invariant complex structure $J_{\bar V}$.

Let $\bar \kv\subset \lg_\IC\times\tilde\kc(\lg)$ be the restriction of the tautological bundle on the Grassmannian to the open subset 
\[\tilde\kc(\lg):= \{\bar V\in \IG(n , \lg_\IC)\mid V\cap \bar V=0\}\]
and consider the manifold \[\tilde\km(\lg):=G\times \tilde\kc(\lg).\]
The group $G$ acts on on the left of $\tilde\km(\lg)$ by  $l_g(h,\bar V)= (gh, \bar V)$ and we can define  the subbundle $\nulleins T\tilde\km(\lg)\subset T\tilde\km(\lg)_\IC$ by 
\[\nulleins T\tilde\km(\lg)\restr{\{g\}\times \tilde\kc(\lg)}:={l_g}_*\bar\kv\oplus \nulleins T\tilde\kc(\lg).\] 

This subbundle gives  an almost complex structure on $\tilde\km(\lg)$ which is integrable over $\kc(\lg)$. So we obtain our desired family by taking the pullback
\[\km(\lg) :=\tilde\km (\lg)\times_{\tilde\kc(\lg)}\kc(\lg).\]

If $\Gamma\subset G$ is a lattice then we can take the quotient of $\km(\lg)$ by the action of $\Gamma$ on the left and we obtain a family $\km(\lg, \Gamma)\to \kc(\lg)$ of compact, complex manifolds such that the fibre over $\bar V\in \kc(\lg)$ is the nilmanifold $M_{\bar V}=( \Gamma\backslash G, J_{\bar V})$. Summarising we have shown the following:

\begin{prop}
Every nilmanifold with left-invariant complex structure $M_{J}$ with fundamental group $\pi(M)\isom \Gamma$ is isomorphic to a nilmanifold in the family $\km(\lg, \Gamma)$.
\end{prop}
\pf We only have to observe that by  \cite{VinGorbShvart}, p.45, Corollary 2.6 the lattice $\pi_1(M)$ determines $\lg$ up to canonical isomorphism, hence $M_{J}$ is biholomorphic to a fibre in the family  $\km(\lg, \Gamma)\to \kc(\lg)$.\qed

There are many natural questions concerning the family $\kc(\lg)$, for example when is it non-empty, smooth, versal and what are the connected components. Catanese and Frediani studied in \cite{catanese04, cat-fred06} the subfamily consisting of principal holomorphic torus bundles over  a torus with fixed dimension of fibre and base, the so called \emph{Appel-Humbert family}, and proved that in some 3-dimensional cases it is a connected component of the Teich\-m\"ul\-ler-Space. The family containing the Iwasawa manifolds was studied by Ketsetzis and Salamon in \cite{ket-sal04}.

\subsection{Kuranishi theory and small deformations}

We will now use deformation theory in the spirit of Kodaira-Spencer \cite{kod-sp58} and Kuranishi \cite{kuranishi62} to study small deformations of left-invariant complex structures on nilmanifolds.

A deformation of a given compact complex manifold $X$ is a flat proper map $\pi:\kx\to \kb$ of (connected) complex spaces, such that all the fibres are smooth manifolds, together with an isomorphism with $X\isom \ky_0=\inverse\pi(0)$ for a point $0\in \kb$. If $\kb$ is  smooth then $\pi$ is just a holomorphic submersion. Kodaira and Spencer showed that first order deformations correspond to elements in $H^1(X, \Theta_X)$ where $\Theta_X$ is the sheaf of holomorphic tangent vectors.

A key result is now the theorem of Kuranishi  which, for a given compact complex manifold $X$, guarantees the existence of a locally complete space of deformations  $\kx\to\Kur(X)$ which is versal at the point corresponding to $X$. In other word, for every deformation $\ky\to \kb$ of $X$ there is a small neighbourhood $\ku$ of $0$ in $\kb$ yielding a diagram
\[\xymatrix{ \ky\restr{\ku}\isom f^*\kx \ar[d]\ar[r] & \kx \ar[d]\\
\ku\ar[r]^f&\Kur(X),}\]
and in addition the differential of $f$ at $0$ is unique.

The Kuranishi family $\Kur(X)$ hence parametrises all sufficiently small deformations of $X$. In general the map $f$ will not be unique which is roughly due to the existence of automorphisms.

In order to study small deformations first of all we need a good description of the cohomology of the tangent bundle.

By a theorem of Nomizu \cite{nomizu54} the de Rham cohomology of  a nilmanifold can be calculated using invariant differential forms and is isomorphic to the cohomology of the complex 
\[0\to \lg^*\overset{d}{\to}\Lambda^2\lg^* \overset{d}{\to}\Lambda^3\lg^* \overset{d}{\to}\dots\]

The question if the Dolbeault-cohomology of compact nilmanifolds with left-invariant complex structure can be calculated using invariant differential forms has been addressed by  Console and Fino in \cite{con-fin01} and  Cordero, Fernandez, Gray and Ugarte in \cite{cfgu00}. We restate their results using the notation from Sections \ref{set-up} and \ref{liedolbeault}:

\begin{theo}\label{citedolbeault}
Let $\Gamma\backslash G=M$ be a real nilmanifold with Lie-algebra $\lg$. Then there is a dense open subset $U$ of the space $\kc(\lg)$ of all left-invariant complex structures on $M$ such that for all $J\in U$  we have an isomorphism 
\[\iota_J:H^{p,q}((\lg,J),\IC)\to H^{p,q}(M_J),\]
on the associated nilmanifold  with left-invariant complex structure $M_J$, 
where we consider $\IC$ as the trivial $\lg_\IC$-module (\cite{con-fin01}, Theorem A).
In addition this holds true in the following cases:
\begin{itemize}
\item  The complex structure $J$ is $\Gamma$-rational. (\cite{con-fin01}, Theorem B).
\item The complex structure $J$ is abelian \cite{con-fin01}.
\item The complex structure $J$ is bi-invariant, $G$ is a complex Lie-group  and $M_J$ is complex parallelisable \cite{sakane76, con-fin01}.
\item The complex manifold $M_J$ has the structure of an iterated principal holomorphic torus bundle \cite{cfgu00}.
\end{itemize}
\end{theo}

The idea of the proof is the following: as long as $M_J$ can be given a structure of iterated bundle with a good control over the cohomology of the base and of the fibre one can use the Borel spectral sequence for Dolbeault-cohomology in order to get an inductive proof. This is the case if the complex structure is $\Gamma$-rational or $M_J$ is an iterated principal holomorphic bundle. This yields the result on a dense subset of the space of invariant complex structures and Console and Fino then show that the property \emph{"The map $\iota_J$ is an isomorphism." } is stable under small deformations.

It is an open question if $\iota_J$ is an isomorphism for every left-invariant complex structure on a nilmanifold.

The work on Lie-algebra Dolbeault-cohomology in Section \ref{LDC} now allows us to compute the cohomology of the holomorphic tangent bundle resp. tangent sheaf:
\begin{cor}
Under the same conditions as in Theorem \ref{citedolbeault} the inclusion
\[\iota:H^{p}_{\delbar}((\lg,J), \lg)\to H^{p}(M_J, \kt_{M_J})\isom H^p(M_J, \Theta_{M_J})\]
is an isomorphism. Here we consider $\lg$ as an integrable $\lg$-module under the adjoint representation.
\end{cor}
\pf This is  Corollary \ref{isoduality} applied to the holomorphic tangent bundle of $M_J$.\qed

The same was proved for 2-step nilmanifolds with abelian complex structure in  \cite{mpps06} and for abelian complex structures in general in \cite{con-fin-poon06}. Hence we can extend the theorem proved there:
\begin{theo}\label{invariantdeformation}
Let $M_J$ be a nilmanifold with left-invariant complex structure of type $(\lg, \Gamma)$ such that 
\[\iota:H^{1,q}((\lg,J),\IC)\to H^{1,q}(M_J)\]
 is an isomorphism for all $q$. Then all small deformations of the complex structure $J$ are again left-invariant complex structures. More precisely, the Kuranishi family contains only left-invariant complex structures. 
\end{theo}

\pf By the work of Kuranishi, the small deformations of $M_J$ are governed by the differential graded algebra $\ka^*_{M_J}(\kt_M)$ of differential forms with values in $\kt_M$. By the above corollary the inclusion $\Lambda^*\nulleins{\lg^*} \tensor \einsnull\lg\subset \ka^*_{M_J}(\kt_M)$ is a quasi-isomorphism and hence induces an isomorphism of corresponding deformation spaces.

We spell this out more in detail following Kuranishi's inductive method on harmonic forms in order to give a description of the Kuranishi space. Note that this has already been done in \cite{mpps06} in the context of abelian complex structures. We choose an invariant, compatible hermitian structure on $M$ as in Section \ref{invariantcohomology}. Recall that the Shouten bracket is defined by 
\begin{gather*}
[\cdot, \cdot]: H^1(M, \kt_M)\times H^1(M, \kt_M) \to H^2(M, \kt_M)\\
[\bar\omega\tensor V, \bar\omega' \tensor V] := \bar \omega' \wedge L_{V'}\bar \omega\tensor V+ \bar\omega \wedge L_{V}\bar\omega'\tensor V'+\bar \omega\wedge \bar \omega ' \tensor [V,V']
\end{gather*}
where $L$ is the Lie derivative, i.e. $L_V\bar\omega'= i_V\circ d \bar\omega'+d\circ i_V\bar\omega'$. By assumption we can represent every class in $H^1(M, \kt_M)$ by an  element in $\nulleins\kh(\lg, \lg)$ which can be considered as an invariant, harmonic differential form on $M$ with respect to the hermitian structure.

Let $G$ be Green's operator which inverts the Laplacian on the orthogonal complement of the harmonic forms. By construction $G$ maps invariant forms to invariant forms since the Laplacian has this property.
Let  $\eta_1, \dots, \eta_m$ be a basis for $\nulleins\kh(\lg, \lg)$ and consider the equation
\[\phi(t)=\sum_{i=1}^m \eta_i t_i +\frac{1}{2} \delbar^* G[\phi(t), \phi(t)].\]

It has a formal power series solution with values in $\nulleins{\lg^*}\tensor \einsnull \lg$ which is given inductively by \[\phi_1(t)=\sum_{i=1}^m \eta_i t_i\text{ and }\phi_r(t)= \frac{1}{2} \sum_{s=1}^{r-1} \delbar^* G[\phi_s(t), \phi_{r-s}(t)].\]
 Note that by construction $\phi(t)$ is left-invariant.

By Kuranishi theory (see e.g. \cite{catanese88}, p. 11) this series converges for small $t$ and there is a complete family of deformations of $M$ over the base
\[B:= \{ t\in B_\epsilon(0)\mid \delbar \phi(t)-\frac{1}{2}[\phi(t), \phi(t)]=0\}.\]
If $\xi_1, \dots, \xi_k$ is a basis of $\kh^{0,2}(\lg, \lg)$ then we can use the inner product $(\cdot, \cdot)$ on $\Lambda^2\nulleins{\lg^*}\tensor \einsnull \lg$ to  describe $B$  as the zero locus of the functions
\[g_i(t)= (\xi_i, [\phi(t), \phi(t) ]),\qquad i=1,\dots ,k.\]
The complex structure over a point $\eta=\sum_{i=1}^m \eta_i t_i\in B$ is determined by
\[ \nulleins{(TM_\eta)}= (id+\phi(t)) \nulleins{TM}.\]
In particular the complex structure is left-invariant since this is true for $\phi(t)$ and $\nulleins{TM}$.\qed

The Kuranishi space has been described more in detail in special cases. If the complex structure is abelian it is often  smooth \cite{mpps06, con-fin-poon06} and if $M_J$ is complex parallelisable it is cut out by polynomial equations but usually singular and reducible \cite{rollenske08a}. We believe that this is true also for general nilmanifolds. 

Beyond small deformations one can look at deformations in the large which have been studied in \cite{rollenske08d}.

% 
% 
% %\bibliographystyle{alpha}
% \bibliographystyle{halpha}
% \bibliography{diss.bib}

\begin{thebibliography}{CFGU00}

\bibitem[And06]{math.DG/0610768}
Adrian Andrada.
\newblock {Complex product structures on 6-dimensional nilpotent Lie algebras},
  2006, arXiv:math.DG/0610768.

\bibitem[BG88]{ben-gor88}
Chal Benson and Carolyn~S. Gordon.
\newblock K\"ahler and symplectic structures on nilmanifolds.
\newblock {\em Topology}, 27(4):513--518, 1988.

\bibitem[Cat88]{catanese88}
F.~Catanese.
\newblock Moduli of algebraic surfaces.
\newblock In {\em Theory of moduli (Montecatini Terme, 1985)}, volume 1337 of
  {\em Lecture Notes in Math.}, pages 1--83. Springer, Berlin, 1988.

\bibitem[Cat04]{catanese04}
Fabrizio Catanese.
\newblock Deformation in the large of some complex manifolds. {I}.
\newblock {\em Ann. Mat. Pura Appl. (4)}, 183(3):261--289, 2004.

\bibitem[CE48]{che-eil48}
Claude Chevalley and Samuel Eilenberg.
\newblock Cohomology theory of {L}ie groups and {L}ie algebras.
\newblock {\em Trans. Amer. Math. Soc.}, 63:85--124, 1948.

\bibitem[CF01]{con-fin01}
S.~Console and A.~Fino.
\newblock Dolbeault cohomology of compact nilmanifolds.
\newblock {\em Transform. Groups}, 6(2):111--124, 2001.

\bibitem[CF06]{cat-fred06}
Fabrizio Catanese and Paola Frediani.
\newblock Deformation in the large of some complex manifolds. {II}.
\newblock In {\em Recent progress on some problems in several complex variables
  and partial differential equations}, volume 400 of {\em Contemp. Math.},
  pages 21--41. Amer. Math. Soc., Providence, RI, 2006.

\bibitem[CFGU00]{cfgu00}
Luis~A. Cordero, Marisa Fern{\'a}ndez, Alfred Gray, and Luis Ugarte.
\newblock Compact nilmanifolds with nilpotent complex structures: {D}olbeault
  cohomology.
\newblock {\em Trans. Amer. Math. Soc.}, 352(12):5405--5433, 2000.

\bibitem[CFP06]{con-fin-poon06}
S.~Console, A.~Fino, and Y.~S. Poon.
\newblock Stability of abelian complex structures.
\newblock {\em Internat. J. Math.}, 17(4):401--416, 2006.

\bibitem[CG90]{Cor-Green}
Lawrence~J. Corwin and Frederick~P. Greenleaf.
\newblock {\em Representations of nilpotent {L}ie groups and their
  applications. {P}art {I}}, volume~18 of {\em Cambridge Studies in Advanced
  Mathematics}.
\newblock Cambridge University Press, Cambridge, 1990.
\newblock Basic theory and examples.

\bibitem[CP07]{0708.3442v2}
Richard Cleyton and Yat~Sun Poon.
\newblock Differential gerstenhaber algebras associated to nilpotent algebras,
  2007, arXiv:0708.3442v2 [math.AG].

\bibitem[Has89]{hasegawa89}
Keizo Hasegawa.
\newblock Minimal models of nilmanifolds.
\newblock {\em Proc. Amer. Math. Soc.}, 106(1):65--71, 1989.

\bibitem[Huy05]{Huybrechts}
Daniel Huybrechts.
\newblock {\em Complex geometry}.
\newblock Universitext. Springer-Verlag, Berlin, 2005.

\bibitem[KN69]{Kob-NumII}
Shoshichi Kobayashi and Katsumi Nomizu.
\newblock {\em Foundations of differential geometry. {V}ol. {II}}.
\newblock Interscience Tracts in Pure and Applied Mathematics, No. 15 Vol. II.
  Interscience Publishers John Wiley \& Sons, Inc., New York-London-Sydney,
  1969.

\bibitem[Kna02]{Knapp}
Anthony~W. Knapp.
\newblock {\em Lie groups beyond an introduction}, volume 140 of {\em Progress
  in Mathematics}.
\newblock Birkh\"auser Boston Inc., Boston, MA, second edition, 2002.

\bibitem[KS58]{kod-sp58}
K.~Kodaira and D.~C. Spencer.
\newblock On deformations of complex analytic structures. {I}, {II}.
\newblock {\em Ann. of Math. (2)}, 67:328--466, 1958.

\bibitem[KS04]{ket-sal04}
Georgios Ketsetzis and Simon Salamon.
\newblock Complex structures on the {I}wasawa manifold.
\newblock {\em Adv. Geom.}, 4(2):165--179, 2004.

\bibitem[Kur62]{kuranishi62}
M.~Kuranishi.
\newblock On the locally complete families of complex analytic structures.
\newblock {\em Ann. of Math. (2)}, 75:536--577, 1962.

\bibitem[MPPS06]{mpps06}
C.~Maclaughlin, H.~Pedersen, Y.~S. Poon, and S.~Salamon.
\newblock Deformation of 2-step nilmanifolds with abelian complex structures.
\newblock {\em J. London Math. Soc. (2)}, 73(1):173--193, 2006.

\bibitem[Nom54]{nomizu54}
Katsumi Nomizu.
\newblock On the cohomology of compact homogeneous spaces of nilpotent {L}ie
  groups.
\newblock {\em Ann. of Math. (2)}, 59:531--538, 1954.

\bibitem[Rol07a]{rollenske07}
S\"onke Rollenske.
\newblock {\em Nilmanifolds: Complex structures, geometry and deformations}.
\newblock PhD thesis, Universit\"at Bayreuth, 2007.

\bibitem[Rol07b]{rollenske07a}
S\"onke Rollenske.
\newblock {The {F}r{\"o}licher spectral sequence can be arbitrarily non
  degenerate}, \newblock {\em Math. Ann.}, 341(3):623--628, 2008, arXiv:0709.0481.

\bibitem[Rol08a]{rollenske08d}
S\"onke Rollenske.
\newblock Geometry of nilmanifolds with left-invariant complex structure and
  deformations in the large, 2008, to appear in Proc. LMS.
\newblock 

\bibitem[Rol08b]{rollenske08a}
S\"onke Rollenske.
\newblock The {K}uranishi space of complex parallelisable nilmanifolds, 2008,
  arXiv:0803.2048v1, to appear in JEMS.

\bibitem[Sak76]{sakane76}
Yusuke Sakane.
\newblock On compact complex parallelisable solvmanifolds.
\newblock {\em Osaka J. Math.}, 13(1):187--212, 1976.

\bibitem[Sal01]{salamon01}
S.~M. Salamon.
\newblock Complex structures on nilpotent {L}ie algebras.
\newblock {\em J. Pure Appl. Algebra}, 157(2-3):311--333, 2001.

\bibitem[VGS00]{VinGorbShvart}
E.~B. Vinberg, V.~V. Gorbatsevich, and O.~V. Shvartsman.
\newblock Discrete subgroups of {L}ie groups [ {MR}0968445 (90c:22036)].
\newblock In {\em Lie groups and Lie algebras, II}, volume~21 of {\em
  Encyclopaedia Math. Sci.}, pages 1--123, 217--223. Springer, Berlin, 2000.

\bibitem[Wei94]{Weibel}
Charles~A. Weibel.
\newblock {\em An introduction to homological algebra}, volume~38 of {\em
  Cambridge Studies in Advanced Mathematics}.
\newblock Cambridge University Press, Cambridge, 1994.

\bibitem[Win98]{winkelmann98}
J{\"o}rg Winkelmann.
\newblock Complex analytic geometry of complex parallelizable manifolds.
\newblock {\em M\'em. Soc. Math. Fr. (N.S.)}, (72-73):x+219, 1998.

\end{thebibliography}
% 

\end{document}